\documentclass[twoside,12pt]{article}
\usepackage{latexsym}
\usepackage{amsmath}
\usepackage{amssymb}
\usepackage{cases}
\usepackage{ascmac}
\usepackage{mathrsfs}
\usepackage{amsthm}

\usepackage{cite}

\usepackage[usenames,dvipsnames]{color}

\title{Nonlinear diffusion equations
\\
with {R}obin boundary conditions 
\\
as asymptotic limits 
of {C}ahn--{H}illiard systems }

\author{Takeshi Fukao\\
Department of Mathematics, Faculty of Education\\
Kyoto University of Education\\
1~Fujinomori, Fukakusa, Fushimi-ku, Kyoto~612-8522 Japan\\
E-mail: \texttt{fukao@kyokyo-u.ac.jp}\\
\and \\ Taishi Motoda\\
Graduate School of Education, 
Kyoto University of Education\\
1~Fujinomori, Fukakusa, Fushimi-ku, Kyoto~612-8522 Japan\\
E-mail: \texttt{motoda.math@gmail.com}}
 
\date{}

\newcommand\testopari{\sc Takeshi Fukao and Taishi Motoda}
\newcommand\testodispari{\sc Nonlinear diffusion equations with {R}obin boundary conditions}
\begin{document}

\maketitle

\begin{abstract}
Condition imposed on the nonlinear terms of a nonlinear diffusion equation with {R}obin boundary condition is the main focus of this paper. 
The degenerate parabolic equations, such as the {S}tefan problem, the {H}ele--{S}haw problem, the porous medium equation and the fast diffusion equation, are included in this class. 
By characterizing this class of equations as an asymptotic limit of the {C}ahn--{H}illiard systems, the growth condition of the nonlinear term can be improved. 
In this paper, the existence and uniqueness of the solution are proved. 
From the physical view point, it is natural that, the {C}ahn--{H}illiard system is treated under the homogeneous {N}eumann boundary condition. 
Therefore, the {C}ahn--{H}illiard system subject to the {R}obin boundary condition looks like pointless. 
However, at some level of approximation, it makes sense to characterize the nonlinear diffusion equations. 

\vspace{2mm}\noindent \textbf{Key words:}~~
{C}ahn--{H}illiard system, degenerate parabolic equation, {R}obin boundary condition, growth condition.

\vspace{2mm}

\end{abstract}

\section{Introduction}
\setcounter{equation}{0}

We consider the initial boundary value problem of a nonlinear diffusion system {\rm (P)}, comprising a parabolic partial differential equation with a {R}obin boundary condition:
\begin{equation*}
	{\rm (P)} \quad 
	\begin{cases}
	\displaystyle \frac{\partial u}{\partial t}-\Delta \xi =g, 
	\quad \xi \in \beta (u) \quad {\rm in~}Q:=\Omega \times (0, T), \\
	\partial_{\boldsymbol{\nu} }\xi +\kappa \xi =h \quad
	{\rm on~}\Sigma :=\Gamma \times (0, T), \\ 
	u(0)=u_{0} \quad {\rm in~}\Omega,
	\end{cases}
\end{equation*}
where $\kappa $ is a positive constant. 
In an asymptotic form, it is characterized as the limit of the {C}ahn--{H}illiard system with a {R}obin boundary condition,
\begin{equation*}
	{\rm (P)}_{\varepsilon }
	\quad 
	\begin{cases}
	\displaystyle \frac{\partial u_\varepsilon }{\partial t}-\Delta \mu_\varepsilon =0 \quad {\rm in~}Q, \\
	\mu_\varepsilon =-\varepsilon \Delta u_\varepsilon 
	+\xi_\varepsilon +\pi_{\varepsilon }(u_\varepsilon )-f, \quad \xi_\varepsilon \in \beta (u_\varepsilon ) 
	\quad {\rm in~}Q, \\
	\partial_{\boldsymbol{\nu} }u_\varepsilon +\kappa u_\varepsilon =0, 
	\quad \partial_{\boldsymbol{\nu} }\mu_\varepsilon +\kappa \mu_\varepsilon =0 
	\quad {\rm on~}\Sigma, 
	\\
	u_\varepsilon (0)=u_{0\varepsilon } \quad {\rm in~}\Omega
	\end{cases}
\end{equation*}
as $\varepsilon \searrow 0$ with $\xi :=\mu +f$, where $0< T< +\infty$, $\Omega $ is a bounded domain of $\mathbb{R}^{d}$ $(d=2,3)$ with smooth boundary $\Gamma :=\partial \Omega $, $\partial_{\boldsymbol{\nu}}$ denotes the outward normal derivative on $\Gamma $, $\Delta $ is the {L}aplacian.
Functions $g: Q\to \mathbb{R}$, $h: \Sigma \to \mathbb{R}$, $u_{0}: \Omega \to \mathbb{R}$ and $u_{0\varepsilon }: \Omega \to \mathbb{R}$ are given as the boundary and initial data.
$f: Q\to \mathbb{R}$ is constructed from $g$ and $h$ later.
Moreover, in the nonlinear diffusion term, $\beta $ is a maximal monotone graph and $\pi_{\varepsilon }$ is an anti-monotone function that tends to 0 as $\varepsilon \searrow 0$. 
It is well known that the {C}ahn--{H}illiard system is characterized by the nonlinear term $\beta +\pi_{\varepsilon }$, a simple example being $\beta (r)=r^{3}$ and $\pi_{\varepsilon }(r)=-\varepsilon r$ for all $r\in \mathbb{R}$. 
In this way, we choose a suitable $\pi_{\varepsilon }$ that depends on the definition of $\beta $ yielding the structure of the {C}ahn--{H}illiard system for ({\rm P})$_\varepsilon $.
Alternatively, by choosing a suitable $\beta $, the degenerate parabolic equation {\rm (P)} characterizes various types of nonlinear problems, such as the {S}tefan problem, {H}ele-{S}haw problem, porous medium equation, and fast diffusion equation (see, e.g., \cite[pp.6935--6937, Examples]{CF16}). 

In analyzing the well-posedness of system {\rm (P)}, there are two standard approaches, the ``$L^1$-approach'' and the `{H}ilbert space approach'' (see, e.g., \cite{Aka09}, \cite[Chapter~5]{Bar10}). 
With respect to the ``{H}ilbert space approach'', the pioneering result \cite{Dam77} concerns the enthalpy formulation for the {S}tefan problem with a {D}irichelet--{R}obin boundary condition, essentially of {R}obin-type. 
Afterwards, the {N}eumann boundary condition was treated in \cite{Ken90, KL05} and the dynamic boundary condition in \cite{Aik93, Aik95, Aik96, FM17}. 
See also \cite{KY17, KY17b, FKY17} for a more general space setting. 
For all of these results, the growth condition for $\beta $ is a very important assumption, such as 
\begin{equation*}
	\widehat{\beta }(r)\geq c_{1}r^{2}-c_{2} \quad {\rm for~all~}r\in \mathbb{R}, 
\end{equation*}
where $c_{1}$ and $c_{2}$ are positive constants, 
and $\widehat{\beta }$ is a proper lower semicontinuous convex function 
satisfying the subdifferential form $\partial_{\mathbb{R}}\widehat{\beta }=\beta $. 
However, it is too restricted in regard to application; 
indeed, for fast diffusion or nonlinear diffusion of {P}enrose--{F}ife type are excluded. 
A drawback of the ``{H}ilbert space approach'' compared with the ``$L^1$-approach'' is detailed in \cite[Chapter~5]{Bar10}. 
With that as motivation, 
the improvement of the growth condition subject to the {R}obin boundary condition, 
was studied in \cite{DK99} using a certain technique called the ``lower semicontinuous convex extension''. 
For porous medium, that is, $\beta (r)=|r|^{q-1}r$, $(q > 1)$, 
a different approach to the doubly nonlinear evolution equation was studied in \cite{Aka09}. 
Regarding a recent result, 
the characterization of the nonlinear diffusion equation as an asymptotic limit of the {C}ahn--{H}illiard system with dynamic boundary conditions was introduced in \cite{Fuk16, Fuk16b}, 
and the same problem subject to the {N}eumann boundary condition was given in \cite{CF16}. 
In these instances, we do not need any growth condition; 
see \cite[Chapter~6]{CF16} or \cite{Fuk16b}. 
This is one of the big advantages of this approach; 
indeed, we do not need techniques such as the lower semicontinuous convex extension by \cite{DK99}.

The main objective of this paper is to improve on the pioneering result give in \cite{Dam77, DK99}, more precisely, an improvement of the growth condition subject to the {R}obin boundary condition without using the lower semicontinuous convex extension. 
Up to a certain level of approximation, we consider the {C}ahn--{H}illiard system subject to the {R}obin boundary condition (see, e.g., \cite{Mil17}). 
From a physical perspective, 
the {C}ahn--{H}illiard system is more naturally treated under the homogeneous {N}eumann boundary condition. 
Therefore the {C}ahn--{H}illiard system with the {R}obin boundary condition imposed looks to be without points. 
However, up to a given level of approximation, characterizing the nonlinear diffusion equations with the {R}obin boundary condition imposed does make sense. 
Moreover, we obtain the order of convergence for the solutions of {\rm (P)} with the {R}obin boundary condition as that for {\rm (P)} with the {N}eumann boundary condition, that is, letting $\kappa $ tend to $0$, where $\kappa $ is the constant in the boundary condition. 

The outline of the paper is as follows. 
In Section~2, the main theorems are stated. 
For this purpose, we present the notation used in this paper and define a suitable duality map and the $H^{1}$-norm equivalent to the standard norm. 
Next, we introduce the definition of a weak solution of {\rm (P)} and {\rm (P)}$_{\varepsilon }$;
the principal theorems are then given. 
In Section~3, to prove the convergence theorem, we deduce the uniform estimates of the approximate solution of {\rm (P)}$_{\varepsilon }$. 
We use {M}oreau--{Y}osida regularization of $\widehat{\beta }$ employing the second-order approximate of parameter $\lambda $. 
In Section~4, to obtain the weak solution of {\rm (P)}$_{\varepsilon }$, we first pass to the limit $\lambda \searrow 0$. 
Second, we prove the existence of weak solutions by passing to the limit $\varepsilon \searrow 0$. 
We also discuss the uniqueness of solutions. 
In Section~5, we improve the assumption for $\beta $ subject to a strong assumption for the heat source $f$. 
From this results, we can avoid the growth condition for $\beta $. 
In Section~6, we obtain the order of convergence related to the {N}eumann problem from the {R}obin problem as $\kappa \searrow 0$. 

Table of contents:
\begin{itemize}
 \item[1.] Introduction
 \item[2.] Main results
\begin{itemize}
 \item[2.1.] Notation
 \item[2.2.] Definition of the solution and main theorem
\end{itemize}
 \item[3.] Approximate problem and uniform estimates
\begin{itemize}
 \item[3.1.] Approximate problem for {\rm (P)}$_{\varepsilon }$
 \item[3.2.] Uniform estimates
\end{itemize}
 \item[4.] Proof of convergence theorem
\begin{itemize}
 \item [4.1.] Passage to the limit $\lambda \searrow 0$
 \item [4.2.] Passage to the limit $\varepsilon \searrow 0$
\end{itemize}
 \item [5.] Improvement of the results
 \item [6.] Asymptotic limits to solutions of {N}eumann problem
 \item [] Appendix
\end{itemize}

\section{Main results}
\setcounter{equation}{0}

In this section, we state the main theorem. Hereafter, $\kappa $ is a positive constant. 

\subsection{Notation}

We employ spaces $H:=L^{2}(\Omega )$ and $H_{\Gamma }:=L^{2}(\Gamma )$, 
with standard norms $|\cdot |_{H}$ and $|\cdot |_{H_{\Gamma }}$, 
along with inner products $(\cdot, \cdot )_{H}$ and $(\cdot, \cdot )_{H_{\Gamma }}$, respectively. 
We also use the space $V:=H^{1}(\Omega )$ with norm $|\cdot |_{V}$ and inner product $(\cdot, \cdot )_{V}$, 
\begin{equation*}
	|z|_{V}:=\bigl(|\nabla z|_{H^{d}}^{2}
	+ \kappa |z|_{H_{\Gamma }}^{2}\bigr)^{\frac{1}{2}}
	\quad {\rm for~all~} z \in V.
\end{equation*}
By virtue of the {P}oincar\'e inequality and the trace theorem, there exist positive constants $c_{\rm P}, c_{\rm P}'$ with $c_{\rm P} < c_{\rm P}'$ such that 
\begin{equation*}
	c_{\rm P} \| z \|_{V}^2 
	\le | z |_{V}^2
	\le c_{\rm P}' \| z \|_{V}^2
	\quad {\rm for~all~} z \in V,
\end{equation*}
where $\| \cdot \|_V$ stands for the standard $H^{1}$-norm. Moreover, we set 
\begin{equation*}
	W:=\bigl\{ z\in H^{2}(\Omega ) \ : \ \partial_{\nu }z+\kappa z=0
	\quad {\rm a.e.~on~}\Gamma \bigr\}. 
\end{equation*}
The symbol $V^{*}$ denotes the dual space of $V$; the duality pairing between $V^{*}$ and $V$ is denoted by $\langle \cdot, \cdot \rangle_{V^{*}, V}$. 
Also, for all $z\in V$, let $F: V\to V^{*}$ be the duality mapping defined by
\begin{equation*}
	\langle Fz, 
	\tilde{z}
	\rangle_{V^{*}, V}
	:=\int_{\Omega }\nabla z\cdot \nabla \tilde{z}dx+\kappa \int_{\Gamma }z\tilde{z}d\Gamma 
	\quad {\rm for~all~}
	z,\tilde{z}\in V.
\end{equation*}
Moreover, we define an inner product of $V^{*}$ by
\begin{equation*}
	(z^{*}, \tilde{z}^{*})_{V^{*}}
	:=\langle z^{*}, F^{-1}\tilde{z}^{*}
	\rangle_{V^{*}, V} 
	\quad 
	{\rm for~all~} 
	z^{*}, \tilde{z}^{*}\in V^{*}.
\end{equation*}
Then, the dense and compact embeddings
$V\mathop{\hookrightarrow} \mathop{\hookrightarrow}
H\mathop{\hookrightarrow} \mathop{\hookrightarrow} V^*$
hold; that is, $(V, H, V^*)$ is a standard {H}ilbert triplet.
As a remark, for all of these settings, $\kappa >0$ is essential. 
For a {N}eumann boundary condition, $\kappa =0$, see \cite{CF16}. 

\subsection{Definition of the solution and main theorem}

We next define our solution for {\rm (P)} and then state the main theorem. 

\paragraph{Definition 2.1.}
{\it The pair $(u, \xi )$ is called the weak solution of {\rm (P)} if 
\begin{gather*}
	u\in H^{1}(0, T;V^{*})\cap L^{\infty }(0, T;H), \quad \xi \in L^{2}(0, T;V), \label{muxi} \\
	\xi \in \beta (u) \quad {\it a.e.~in~}Q, \label{2.20}
\end{gather*}
and they satisfy}
\begin{gather}
	\bigl\langle u'(t), z\bigr\rangle_{V^{*}, V}
	+\int_{\Omega }\nabla \xi (t)\cdot \nabla zdx
	+\kappa \int_{\Gamma }\xi(t) zd\Gamma =
	\int_{\Omega }^{} g(t)z dx +\int_{\Gamma }^{} h(t)z d\Gamma \nonumber \\
	 \quad {\it for~all~}z\in V, \quad 
	 {\it for~a.a.\ } t \in (0,T), \label{2.19'} \\
	u(0)=u_{0} \quad {\it a.e.~in}~\Omega. \label{ic}
\end{gather}

In this definition, the {R}obin boundary condition for $\xi$ is hidden in the weak formulation \eqref{2.19'}.
The strategy behind the proof of the main theorem is the characterization of our nonlinear diffusion equation {\rm (P)} as an asymptotic limit of the {C}ahn--{H}illiard system. 
Therefore, for each $\varepsilon \in (0,1]$, we define the approximate problem of {C}ahn--{H}illiard type with {R}obin boundary condition as follows:

\paragraph{Definition 2.2.}
{\it The triplet $(u_{\varepsilon }, \mu_{\varepsilon }, \xi_{\varepsilon })$ is called the weak solution of {\rm (P)}$_{\varepsilon }$ if 
\begin{gather}
	u_{\varepsilon }\in H^{1}(0, T;V^{*})\cap L^{\infty }(0, T;V) \cap L^2(0,T;W), \label{3.1}\\
	\mu_{\varepsilon }\in L^{2}(0, T;V), \quad \xi_{\varepsilon }\in L^{2}(0, T;H), \nonumber \\
	\xi_{\varepsilon }\in \beta (u_{\varepsilon }) \quad {\it a.e.~in~}Q \nonumber 
\end{gather}
and they satisfy
\begin{gather}
	\bigl\langle u_{\varepsilon }'(t), z\bigr\rangle_{V^{*}, V}
	+\int_{\Omega }\nabla \mu_{\varepsilon }(t)\cdot \nabla zdx
	+\kappa \int_{\Gamma }\mu_{\varepsilon } (t)zd\Gamma =0 \quad {\it for~all~}z\in V, \label{2.13}\\
	\mu_{\varepsilon }(t)
	=-\varepsilon \Delta u_{\varepsilon }(t)
	+\xi_{\varepsilon }(t)+\pi_{\varepsilon } \bigl( u_{\varepsilon }(t) \bigr)
	-f(t) \quad {\it in~}H, \label{2.14}
\end{gather}
for a.a.\ $t\in (0, T)$, with}
\begin{equation*}
	u_\varepsilon (0)=u_{0\varepsilon } \quad {\it a.e.~in}~\Omega.
\end{equation*}

The {R}obin boundary condition for $u_\varepsilon $ is stated with regard to the class of $W$, that is, the regularity \eqref{3.1}; 
that for $\mu_\varepsilon $ is hidden in the weak formulation \eqref{2.13}. 
The {C}ahn--{H}illiard structure is characterized by the nonlinear term $\beta+\pi_\varepsilon $. 
The conditions for these terms are given as assumptions: 
\begin{enumerate}
 \item[(A1)] 
	$\beta :\mathbb{R}\to 2^{\mathbb{R}}$ is a maximal monotone graph, which is
	the subdifferential $\beta =\partial_{\mathbb{R}}\widehat{\beta }$ 
	of some proper lower semicontinuous convex function 
	$\widehat{\beta }: \mathbb{R}\to [0, +\infty ]$ 
	satisfying $\widehat{\beta }(0)=0$ with the effective domain $D(\beta ):=\{ r \in \mathbb{R} : \beta (r) \ne \emptyset \}$; 
 \item[(A2)]
	there exist positive constants $c_{1}$, $c_{2}$ such that 
	$\widehat{\beta }(r)\geq c_{1}r^{2}-c_{2}$ for all $r\in \mathbb{R}$; 
 \item[(A3)]
	$\pi_{\varepsilon }: \mathbb{R}\to \mathbb{R}$ is a Lipschitz continuous function for all $\varepsilon \in (0, 1]$. 
	Moreover, there exist a positive constant $c_{3}$ and strictly increasing continuous function $\sigma : [0, 1]\to [0, 1]$ such that $\sigma (0)=0$, $\sigma (1)=1$, and 
\begin{equation}
	\bigl|\pi_{\varepsilon }(0)\bigr|+
	|\pi_{\varepsilon }'|_{L^{\infty }(\mathbb{R})}\leq c_{3}\sigma (\varepsilon ) 
	\quad {\rm for~all~}\varepsilon \in (0, 1]. \label{2.8}
\end{equation}
\end{enumerate}
In particular, {\rm (A1)} yields $0\in \beta (0)$. 
Assumption {\rm (A2)} is improved in Section~6. 
The assumptions pertaining to the given data are as follows:
\begin{enumerate}
 \item[(A4)]
	$g\in L^{2}(0, T; H)$, $h\in L^{2}(0, T; H_{\Gamma })$;
 \item[(A5)]
	$u_{0}\in H$ with $\widehat{\beta }(u_{0})\in L^{1}(\Omega )$. Moreover, let 
	$u_{0\varepsilon }\in V$; then there exists a positive constant 
	$c_{4}$ such that, for all $\varepsilon \in (0, 1]$, 
\begin{equation}\label{2.9}
	|u_{0\varepsilon }|_{H}^{2}\leq c_{4}, 
	\quad \int_{\Omega }\widehat{\beta }(u_{0\varepsilon })dx\leq c_{4}, 
	\quad \varepsilon |\nabla u_{0\varepsilon }|_{H^{d}}^{2}\leq c_{4}, 
	\quad \varepsilon |u_{0\varepsilon }|_{H_{\Gamma }}^{2}\leq c_{4}. 
\end{equation}
In addition, $u_{0\varepsilon }\to u$ strongly in $H$ as $\varepsilon \searrow 0$ (cf.\ \cite[Lemma~A.1]{CF16}). 
\end{enumerate}

From assumption {\rm (A4)}, we see from the {L}ax--{M}ilgram theorem that there exists a unique function $f\in L^{2}(0, T; V)$ such that
\begin{equation*}
	\int_{\Omega }\nabla f(t)\cdot \nabla zdx
	+\kappa \int_{\Gamma }f(t)zd\Gamma 
	=\int_{\Omega }g(t)zdx+\int_{\Gamma }h(t)zd\Gamma
\end{equation*}
for all $z\in V$ and for a.a.\ $t\in (0, T)$. 
Therefore, introducing new variable $\mu :=\xi -f \in L^2(0,T;V)$, we rewrite the weak formulation \eqref{2.19'} as follows: 
\begin{equation}
	\bigl\langle u'(t), z\bigr\rangle_{V^{*}, V}
	+\int_{\Omega }\nabla \mu (t)\cdot \nabla zdx
	+\kappa \int_{\Gamma }\mu(t) zd\Gamma =0 
	 \quad {\rm for~all~}z\in V, \label{2.19} 
\end{equation}
for a.a.\ $t \in (0,T)$. 

The proof of main theorem follows that in \cite{Fuk16, CF16} for the {R}obin boundary condition. 
The characterization of the nonlinear diffusion equation from the asymptotic limits of 
{C}ahn--{H}illiard system \cite{Fuk16, Fuk16b, CF16, FKY17} furnishes a big advantage 
in regard to the growth condition for $\beta $. 
Because of this, we can improve the result \cite{Dam77} and widen the setting for 
$\beta $ using the different approach described in \cite{DK99} starting from the lower semicontinuous convex extension. 
To do so, we replace assumption {\rm (A4)} by the following {\rm (A6)}: 

\begin{enumerate}
 \item [(A6)]
	$g\in L^{2}(0, T; H)$, $h=0$ a.e.\ on $\Sigma $. 
\end{enumerate}

Then, we see that there exists a unique function $f\in L^{2}(0, T; V)$ such that,
\begin{equation*}
	\int_{\Omega }\nabla f(t)\cdot \nabla zdx
	+\kappa \int_{\Gamma }f(t)zd\Gamma =\int_{\Omega }g(t)zdx
	\quad {\rm for~all~} z \in V,
\end{equation*}
for a.a.\ $t\in (0, T)$. 
Now, taking a test function $z \in \mathcal{D}(\Omega )$, we have $-\Delta f(t)=g(t)$ in $\mathcal{D}'(\Omega )$. Specifically, by comparison, $-\Delta f(t)=g(t)$ in $H$. 
This yields 
\begin{equation*}
	\partial_{\boldsymbol{\nu }}f(t)+\kappa f(t)=0 \quad {\rm a.e.~on~}\Gamma, 
\end{equation*}
for a.a.\ $t\in (0, T)$. 
Therefore, under assumption {\rm (A6)}, we have $-\Delta f\in L^{2}(0, T; H)$ and $\partial_{\boldsymbol{\nu }}f\in L^{2}(0, T; H_{\Gamma })$. 
These higher regularities are essential to improve the growth condition {\rm (A2)}. 

The well-posedness of {C}ahn--{H}illiard system has been treated in many studies (see, e.g., \cite{Mil17}). 
In regard to the abstract theorem of the evolution equation, we refer the reader to \cite{KNP95, Kub12}. 
Based on these results, we obtain the following proposition: 

\paragraph{Proposition 2.1.}{\it Given assumptions {\rm (A1)}--{\rm (A5)} or {\rm (A1)}, {\rm (A3)} with $\sigma (\varepsilon )=\varepsilon ^{1/2}$, {\rm (A5)}, {\rm (A6)}, then for each $\varepsilon \in (0, 1]$, there exists a unique weak solution of {\rm (P)}$_\varepsilon $. } \\

This proposition implies that, for the well-posedness of the {C}ahn--{H}illiard system, the growth condition {\rm (A2)} is not essential. It can be recovered by the strong assumption {\rm (A6)}.
The proof of this proposition is given in Section~4. Indeed, we can prove this proposition by considering the approximate problem given in Proposition~3.1.

Our main theorem is now given. 

\paragraph{Theorem 2.1.}
{\it With assumptions {\rm (A1)}--{\rm (A5)}, for each $\varepsilon \in (0, 1]$, let $(u_{\varepsilon }, \mu_{\varepsilon }, \xi_{\varepsilon })$ be the weak solution of {\rm (P)}$_\varepsilon $ obtained in Proposition~2.1. 
Then, there exists a weak solution $(u,\xi )$ of {\rm (P)} and $(u, \xi)$ characterized by $(u_{\varepsilon }, \mu_{\varepsilon }, \xi_{\varepsilon })$ in the following sense: 
\begin{gather*}
	u_\varepsilon \to u \quad {\it strongly~in~} C\bigl( [0,T];V^* \bigr) 
	\quad {\it and~weakly~star~in~}H^1(0,T;V^*) \cap L^\infty (0,T;H),
	\\
	\xi_\varepsilon \to \xi
	\quad {\it weakly~in~}L^2(0,T;H), \\
	\mu_\varepsilon \to \mu:=\xi -f 
	\quad {\it weakly~in~}L^2(0,T;V)
\end{gather*}
as $\varepsilon \searrow 0$. 
Moreover, the component $u$ of the solution of {\rm (P)} is uniquely determined. Also, if $\beta $ is single-valued, then the component $\xi $ of the solution of {\rm (P)} is also unique. 
}\\

The second theorem relates to improving the well-posedness result of \cite{Dam77, DK99}. 

\paragraph{Theorem 2.2.}
{\it Given assumptions {\rm (A1)}, {\rm (A3)} with $\sigma (\varepsilon )=\varepsilon ^{1/2}$, {\rm (A5)}, {\rm (A6)}, the same statement as in Theorem~2.1 holds.
} 

\section{Approximate problem and uniform estimates}
\setcounter{equation}{0}

The proof of the main theorems exploits the characterization of the nonlinear diffusion equation through the asymptotic limits of the {C}ahn--{H}illiard system \cite{Fuk16, CF16}.
To apply it, we consider the second-order approximation of the nonlinear term $\beta $. 
At this level of approximation, we obtain a uniform estimate independent of the approximation parameters. 

\subsection{Approximate problem for (P)$_{\varepsilon }$}

We consider an approximate problem to show the well-posedness of (P)$_{\varepsilon }$. 
For each $\lambda \in (0,1]$, we define $\beta_{\lambda }: \mathbb{R}\to \mathbb{R}$ by 
\begin{equation*}
	\beta_{\lambda }(r)
	:=\frac{1}{\lambda }\bigl(r-J_{\lambda }(r)\bigr) \quad {\rm for~all~}r\in \mathbb{R}, 
\end{equation*}
where the resolvent operator $J_{\lambda }: \mathbb{R}\to \mathbb{R}$ is given by 
\begin{equation*}
	J_{\lambda }(r):=(I+\lambda \beta )^{-1}(r) \quad {\rm for~all~}r\in \mathbb{R}. 
\end{equation*}
Also, we define the {M}oreau--{Y}osida regularization $\widehat{\beta }_{\lambda }$ of $\widehat{\beta }: \mathbb{R}\to \mathbb{R}$ by
\begin{equation*}
	\widehat{\beta }_{\lambda }(r)
	:=\inf_{s\in \mathbb{R}}\left\{\frac{1}{2\lambda }|r-s|^{2}+\widehat{\beta }(s)\right\}
	=\frac{1}{2\lambda }\bigl| r-J_{\lambda }(r) \bigr|^{2}
	+ \widehat{\beta } \bigl( J_{\lambda }(r) \bigr) 
	\quad {\rm for~all~} r\in \mathbb{R}.
\end{equation*}
Now, we consider the problem (P)$_{\varepsilon, \lambda }$ for the viscous {C}ahn--{H}illiard like system:
\begin{equation*}
	{\rm (P)}_{\varepsilon, \lambda }
	\quad 
	\begin{cases}
	\displaystyle \frac{\partial u_{\varepsilon, \lambda }}{\partial t}
	-\Delta \mu_{\varepsilon, \lambda }=0 \quad {\rm a.e.\ in~}Q, \\
	\displaystyle \mu_{\varepsilon, \lambda }
	=\lambda \frac{\partial u_{\varepsilon, \lambda }}{\partial t}
	-\varepsilon \Delta u_{\varepsilon, \lambda }
	+\beta_{\lambda }(u_{\varepsilon, \lambda })
	+\pi_{\varepsilon }(u_{\varepsilon, \lambda })-f \quad {\rm a.e.\ in~}Q, \\
	\partial_{\boldsymbol{\nu} }u_{\varepsilon, \lambda }
	+\kappa u_{\varepsilon, \lambda }
	=0, \quad \partial_{\boldsymbol{\nu} }\mu_{\varepsilon, \lambda }
	+\kappa \mu_{\varepsilon, \lambda }=0 \quad {\rm a.e.\ on~}\Sigma,\\
	u_{\varepsilon, \lambda }(0)=u_{0\varepsilon } \quad {\rm a.e.\ in~}\Omega. 
	\end{cases}
\end{equation*}
Define $A:D(A) \to H$ by $A u =-\Delta u$ in $H$ with $D(A)=W$; the treatment of $A$ is given in Appendix. 
From the well-known abstract theory of the doubly nonlinear evolution equation \cite{CV90}, 
we obtain the following well-posedness result (see, also \cite{CF15, KNP95, Kub12}): 

\paragraph{Proposition 3.1.}{\it 
For each $\lambda \in (0, 1]$, there exists a unique 
\begin{equation*}
	u_{\varepsilon, \lambda }\in H^{1}(0, T; H)\cap L^{\infty }(0, T; V)\cap L^{2}(0, T; W)
\end{equation*}
such that $u_{\varepsilon, \lambda }$ satisfies the following {C}auchy problem}
\begin{align}
	& (\lambda I+F^{-1})u'_{\varepsilon, \lambda }(t)
	+ \varepsilon A
	u_{\varepsilon, \lambda }(t) \nonumber \\
	& \quad =-\beta_{\lambda }\bigl( u_{\varepsilon, \lambda }(t) \bigr)
	-\pi_{\varepsilon } \bigl( u_{\varepsilon, \lambda }(t) \bigr)
	+f(t) \quad {\it in~}H, \quad {\it for~a.a.~} t\in (0, T), \label{1} \\
	& u_{\varepsilon, \lambda }(0)=u_{0\varepsilon } \quad {\it in~}H. \nonumber 
\end{align} 

With this level of abstractness, the {C}ahn--{H}illiard system with {R}obin boundary condition is essentially the same as in previous studies. Therefore, we omit the proof of this proposition. \\

Now, for each $\lambda \in (0, 1]$, we put 
\begin{equation}
	\mu_{\varepsilon, \lambda }(t)
	:=\lambda u_{\varepsilon, \lambda }'(s)
	+\varepsilon A u_{\varepsilon, \lambda }(t) 
	+\beta_{\lambda } \bigl( u_{\varepsilon, \lambda }(t) \bigr)
	+\pi_{\varepsilon } \bigl( u_{\varepsilon, \lambda }(t) \bigr)
	-f(t) \quad {\rm in}~H,
	\label{9}
\end{equation}
for a.a.\ $t\in (0, T)$.
Then, we can rewrite the evolution equation \eqref{1} as 
\begin{equation}
	F^{-1}u'_{\varepsilon, \lambda }(t)+\mu_{\varepsilon, \lambda }(t)=0 \quad {\rm in~}V,
	\label{10}
\end{equation}
for a.a.\ $t\in (0, T)$. 
We remark here that we do not need the projection to $\mu_{\varepsilon, \lambda}$ because of the boundary condition (cf.\ \cite{CF15, KNP95, Kub12}). 
This is different from that for the {N}eumann boundary condition. 

\subsection{Uniform estimates}

To prove the convergence theorem, we now obtain the uniform estimates independent of $\varepsilon, \lambda $. 

\paragraph{Lemma 3.1.}
{\it There exists a positive constant $M_{1}$ and two values $\bar{\lambda }, \bar{\varepsilon }\in (0, 1]$, depending only on the data, such that }
\begin{gather}
	\int_{0}^{t}
	\bigl|u'_{\varepsilon, \lambda }(s)
	\bigr|_{V^{*}}^{2}ds
	+2\lambda \int_{0}^{t}
	\bigl|u'_{\varepsilon, \lambda }(s)\bigr|_{H}^{2}ds
	+\varepsilon \bigl| u_{\varepsilon, \lambda }(t) \bigr|_{V}^{2} \nonumber \\
	{}+\bigl| \widehat{\beta }_{\lambda } \bigl( u_{\varepsilon, \lambda }(t) \bigr) 
	\bigr|_{L^{1}(\Omega )}
	+ \frac{c_{1}}{4} \bigl| u_{\varepsilon, \lambda }(t) \bigr |_{H}^{2} \leq M_{1}, \label{m1} \\
	\int_{0}^{t}\bigl| \mu_{\varepsilon, \lambda }(s) \bigr|_{V}^{2}ds\leq M_{1} \label{m2}
\end{gather}
{\it for all $t\in [0, T], \lambda \in (0, \bar{\lambda }]$ and $\varepsilon \in (0, \bar{\varepsilon }]$. }

\paragraph{Proof}
We test \eqref{1} at time $s\in (0, T)$ for $u_{\varepsilon, \lambda }'(s)\in H$. 
Then, we see that 
\begin{align}
	& \lambda \bigl|u_{\varepsilon, \lambda }'(s)
	\bigr|_{H}^{2}+\bigl|u_{\varepsilon, \lambda }'(s)
	\bigr|_{V^{*}}^{2}
	- \varepsilon
	\bigl( \Delta u_{\varepsilon, \lambda }(s), u_{\varepsilon, \lambda }'(s)\bigr)_{H} 
	\nonumber \\
	& \quad {} 
	+\frac{d}{ds} \int_{\Omega } \widehat{\beta }_{\lambda }
	\bigl( u_{\varepsilon, \lambda }(s) \bigr)dx
	=-\bigl(\pi_{\varepsilon } \bigl( u_{\varepsilon, \lambda }(s) \bigr), u_{\varepsilon, \lambda }'(s)\bigr)_{H}
	+\bigl(f(s), u_{\varepsilon, \lambda }'(s)\bigr)_{H} \label{m1-1}
\end{align}
for a.a.\ $s\in (0, T)$. 
We have now from the boundary condition of $u_{\varepsilon, \lambda }(s)$
\begin{align*}
	- \bigl(
	\Delta u_{\varepsilon, \lambda }(s), u_{\varepsilon, \lambda }'(s)
	\bigr)_{H}
	&= - \int_{\Gamma }\partial_{\boldsymbol{\nu }}u_{\varepsilon, \lambda }(s)u_{\varepsilon, \lambda }'(s)d\Gamma 
	 + \int_{\Omega }\nabla u_{\varepsilon, \lambda }(s)\cdot \nabla u_{\varepsilon, \lambda }'(s)dx \\
	&= \kappa \int_{\Gamma }u_{\varepsilon, \lambda }(s)u_{\varepsilon, \lambda }'(s)d\Gamma 
	+\int_{\Omega }\nabla u_{\varepsilon, \lambda }(s)\cdot \nabla u_{\varepsilon, \lambda }'(s)dx \\
	&= \frac{\kappa }{2}
	\frac{d}{ds}\int_{\Gamma } \bigl|u_{\varepsilon, \lambda }(s) \bigr|^{2}d\Gamma 
	+\frac{1}{2}
	\frac{d}{ds}\int_{\Omega } \bigl|\nabla u_{\varepsilon, \lambda }(s) \bigr|^{2}dx
\end{align*}
for a.a.\ $s\in (0, T)$. 
Integrating \eqref{m1-1}with respect to $s$ over interval $[0, t]$, and using the above, we infer that 
\begin{align}
	&\int_{0}^{t}
	\bigl|u_{\varepsilon, \lambda }'(s)\bigr|_{V^{*}}^{2}
	ds
	+\lambda \int_{0}^{t}
	\bigl|u_{\varepsilon, \lambda }'(s)\bigr|_{H}^{2}
	ds
	+\frac{\varepsilon }{2}
	\bigl|u_{\varepsilon, \lambda }(t)
	\bigr|_{V}^{2}+\int_{\Omega }\widehat{\beta }_{\lambda }
	\bigl( u_{\varepsilon, \lambda }(t) \bigr)dx \nonumber \\
	&=
	\frac{\varepsilon }{2}|u_{0\varepsilon }|_{V}^{2}
	+\int_{\Omega }\widehat{\beta }_{\lambda }(u_{0\varepsilon })dx
	+\int_{\Omega }\widehat{\pi }_{\varepsilon }(u_{0\varepsilon })dx
	-\int_{\Omega }\widehat{\pi }_{\varepsilon }
	\bigl(u_{\varepsilon, \lambda }(t) \bigr)dx \nonumber \\
	&\quad {}+\int_{0}^{t}\bigl\langle u_{\varepsilon, \lambda }'(s), f(s)
	\bigr\rangle_{V^{*}, V}ds \label{3.6}
\end{align}
for all $t\in [0, T]$, where $\widehat{\pi }_{\varepsilon }$ is the primitive of $\pi_{\varepsilon }$ given by $\widehat{\pi }_{\varepsilon }(r):=\int_{0}^{r} \pi_{\varepsilon }(\tau)d\tau \ {\rm for~all~} r\in \mathbb{R}$. 

Now, aided by assumption {\rm (A2)}, we have
\begin{align*}
	\widehat{\beta }_{\lambda }(r)
	& = \frac{1}{2\lambda }
	\bigl|r-J_{\lambda }(r)
	\bigr|^{2}
	+\widehat{\beta }_{\lambda }\bigl(J_{\lambda }(r)\bigr) \\
	& \geq 
	\frac{1}{2\lambda }
	\bigl|r-J_{\lambda }(r)\bigr|^{2}+c_{1}
	\bigl|J_{\lambda }(r)\bigr|^{2}-c_{2} \quad {\rm for~all~}r\in \mathbb{R}. 
\end{align*}
Hence, putting $\bar{\lambda }:=\min \{1, 1/(2c_{1})\}$, we see that $c_{1}/2\leq 1/(4\bar{\lambda })$. 
It follows that 
\begin{align}
	\int_{\Omega }
	\widehat{\beta }_{\lambda }
	\bigl( u_{\varepsilon, \lambda }(s) \bigr)
	&\geq 
	\frac{1}{2}
	\int_{\Omega }\widehat{\beta }_{\lambda } \bigl( u_{\varepsilon, \lambda }(s) \bigr)dx
	+\frac{1}{4\bar{\lambda }}
	\bigl|u_{\varepsilon, \lambda }(s)-J_{\lambda }
	\bigl( 
	u_{\varepsilon, \lambda }(s)
	\bigr)
	\bigr|_{H}^{2} \nonumber \\
	&\quad {} +\frac{c_{1}}{2}
	\bigl|
	J_{\lambda } \bigl( 
	u_{\varepsilon, \lambda }(s) \bigr) \bigr|_{H}^{2}-\frac{c_{2}}{2}|\Omega | \nonumber \\
	&\geq \frac{1}{2}
	\int_{\Omega }\widehat{\beta }_{\lambda }
	\bigl(
	u_{\varepsilon, \lambda }(s) \bigr) dx+\frac{c_{1}}{4}\bigl|u_{\varepsilon, \lambda }(s)\bigr|_{H}^{2}-\frac{c_{2}}{2}|\Omega | \label{3.7}
\end{align}
for a.a.\ $s \in (0, T)$. 
Next, we use the {M}aclaurin expansion and {\rm (A3)} to obtain 
\begin{equation*}
	\bigl| \widehat{\pi }_{\varepsilon }(r) \bigr|
	\leq 
	\bigl| \pi_{\varepsilon }(0) \bigr| 
	|r| 
	+ \frac{1}{2} 
	|\pi_{\varepsilon }'|_{L^{\infty }(\mathbb{R})}
	r^{2}
	\leq c_{3} \sigma (\varepsilon )(1+r^{2})
\end{equation*}
for all $r\in \mathbb{R}$. 
Also, with assumption {\rm (A3)}, there exists $\bar{\varepsilon }\in (0, 1]$ such that 
\begin{equation*}
	\sigma (\varepsilon )\leq \frac{c_{1}}{8c_{3}(1+|\Omega |)}
\end{equation*}
for all $\varepsilon \in (0, \bar{\varepsilon }]$; that is, we deduce that 
\begin{align}
	- \int_{\Omega }
	\widehat{\pi }_{\varepsilon }
	\bigl( u_{\varepsilon, \lambda }(t) \bigr)dx
	&\leq c_{3}\sigma (\varepsilon )\int_{\Omega }
	\bigl(1+\bigl|u_{\varepsilon, \lambda }(t)\bigr|^{2}\bigr)dx \nonumber \\
	&\leq \frac{c_{1}}{8}\bigl(1+\bigl|u_{\varepsilon, \lambda }(t)\bigr|_{H}^{2}\bigr)
	\label{s6}
\end{align}
for all $t \in [0, T]$. 
Hence, applying the {Y}oung inequality to \eqref{3.6} and collecting all of the above, we see that there exists a positive constant $M_{1}$ depending only on $c_{1}$, $c_{2}$, $c_{4}$, $|\Omega |$, and $|f|_{L^{2}(0, T; V)}$, independent of $\varepsilon \in (0, \bar{\varepsilon }], \lambda \in (0, \bar{\lambda }]$ such that \eqref{m1} holds. Finally, we have from \eqref{10}
\begin{equation*}
	\int_{0}^{t} \bigl| \mu_{\varepsilon, \lambda }(s) \bigr|_{V}^{2}ds
	= \int_{0}^{t} \bigl| u'_{\varepsilon, \lambda }(s) \bigr|_{V^*}^{2}\leq M_{1}
\end{equation*}
for all $t \in [0,T]$. \hfill $\Box $

\paragraph{Lemma 3.2.}
{\it There exist positive constants $M_{2}$ and $M_{3}$, independent of $\varepsilon \in (0, \bar{\varepsilon }]$ and $\lambda \in (0, \bar{\lambda }]$, such that }
\begin{gather}
	\int_{0}^{t} \bigl|\beta_{\lambda }
	\bigl( u_{\varepsilon, \lambda }(s) \bigr) \bigr|_{H}^2ds \leq M_{2}, 
	\label{4.23}
	\\
	\int_{0}^{t}
	\bigl|\varepsilon \Delta u_{\varepsilon, \lambda }(s)
	\bigr|_{H}^{2}ds\leq M_{3}
	\label{lemma33b}
\end{gather}
{\it for all} $t\in [0, T]$. 

\paragraph{Proof}
We test \eqref{9} for times $s\in (0, T)$ by $\beta_{\lambda }(u_{\varepsilon, \lambda }(s))\in H $. Then, we see that 
\begin{gather*}
	\bigl( \mu_{\varepsilon, \lambda }(s), 
	\beta_{\lambda }\bigl( u_{\varepsilon, \lambda }(s) \bigr) 
	\bigr)_{H}
	= - \varepsilon \bigl( \Delta u_{\varepsilon, \lambda }(s), 
	\beta_{\lambda }\bigl(u_{\varepsilon, \lambda }( s )\bigr)
	\bigr)_{H} 
	+ \bigl| \beta_{\lambda } \bigl( u_{\varepsilon, \lambda }(s) \bigr)
	\bigr|_{H}^2 \\
	{} + \bigl(\lambda u_{\varepsilon, \lambda }'(s)
	+ \pi_{\varepsilon } \bigl( 
	u_{\varepsilon, \lambda }(s)
	\bigr) - f(s), 
	\beta_{\lambda } \bigl( u_{\varepsilon, \lambda }(s) \bigr)
	\bigr)_{H}
\end{gather*}
for a.a.\ $s\in (0, T)$. 
Calculating the first term on the right-hand side of the above equation, we infer from the boundary condition of $u_{\varepsilon, \lambda }(s)$ that 
\begin{align*}
	-\varepsilon \bigl(\Delta u_{\varepsilon, \lambda }(s), 
	\beta_{\lambda }
	\bigl( u_{\varepsilon, \lambda }(s) \bigr)\bigr)
	& = 
	\varepsilon \int_{\Omega }
	\beta_{\lambda }' \bigl( u_{\varepsilon, \lambda }(s) \bigr)
	\bigl|\nabla u_{\varepsilon, \lambda }(s)\bigr|^{2}dx
	+\kappa \varepsilon 
	\int_{\Gamma }
	u_{\varepsilon, \lambda }(s)
	\beta_{\lambda }\bigl( u_{\varepsilon, \lambda }(s) \bigr)d\Gamma
	\\
	& \ge 0
\end{align*}
for a.a.\ $s\in (0, T)$ because $\beta$ is monotonic. Hence, applying the {Y}oung inequality, we obtain
\begin{align*}
	&\bigl|\beta_{\lambda } \bigl( u_{\varepsilon, \lambda }(s) \bigr)\bigr|_{H}^{2} \\
	&\quad \leq \bigl( \mu_{\varepsilon, \lambda }(s)
	- \lambda u_{\varepsilon, \lambda }'(s)
	- \pi_{\varepsilon } \bigl( u_{\varepsilon, \lambda }(s) \bigr) 
	+ f(s), \beta_{\lambda }\bigl( u_{\varepsilon, \lambda }(s) \bigr) \bigr)_{H} \\
	&\quad \leq \frac{1}{2}
	\bigl|\beta_{\lambda } \bigl( u_{\varepsilon, \lambda }(s) \bigr)
	\bigr|_{H}^{2} 
	+ 2\left( 
	\bigl| \mu_{\varepsilon, \lambda }(s) \bigr|_{H}^{2}
	+\lambda ^{2}\bigl| u_{\varepsilon, \lambda }'(s) \bigr|_{H}^{2}
	+\bigl|\pi_{\varepsilon } \bigl( u_{\varepsilon, \lambda }(s) \bigr)
	\bigr|_{H}^{2}
	+\bigl| f(s) \bigr|_{H}^{2} \right)
\end{align*}
for a.a.\ $s\in (0, T)$. 
Integrating with respect to $s$ over $[0, t]$ and applying Lemma~3.2, we see that there exists a positive constant $M_{2}$ depending only on $M_{1}$, $c_3$, and $|f|_{L^{2}(0, T; H)}$ such that estimate \eqref{4.23} holds. 

Next, by comparing with \eqref{9}, we obtain \eqref{lemma33b} for some positive constant $M_{3}$.
$\hfill $ $\Box $


\section{Proof of convergence theorem}

We next show the existence of solution of (P)$_{\varepsilon }$ by passing to the limit for the approximate problem (P)$_{\varepsilon, \lambda }$ . 

\subsection{Passage to the limit $\lambda \searrow 0$}

\paragraph{Proof of Proposition 2.1.}
\indent
From previous estimates established in Lemmas~3.1 and 3.2, 
we see that there exists a subsequence $\{\lambda_{k}\}_{k\in \mathbb{N}}$ with $\lambda_{k}\searrow 0$ 
as $k\nearrow +\infty $ and some limit functions 
$u_{\varepsilon }\in H^{1}(0, T; V^{*})\cap L^{\infty }(0, T; V)\cap L^2(0,T;W)$, $\mu_{\varepsilon }\in L^{2}(0, T; V)$, and $\xi_{\varepsilon }\in L^{2}(0, T; H)$ such that 
\begin{gather}
	\label{4.27}
	u_{\varepsilon, \lambda_{k}}\to u_{\varepsilon } 
	\quad {\rm weakly~star~in~} H^{1}(0, T; V^{*})\cap L^{\infty }(0, T; V), 
	\\ \label{4.28}
	\lambda_{k}u_{\varepsilon, \lambda_{k}}\to 0 
	\quad {\rm strongly~in~} H^{1}(0, T; H), 
	\\
	\label{4.29}
	\mu_{\varepsilon, \lambda_{k}}\to \mu_{\varepsilon } 
	\quad {\rm weakly~in~} L^{2}(0, T; V), 
	\\
	\label{4.30}
	\beta_{\lambda_{k}}(u_{\varepsilon, \lambda_{k}})\to \xi_{\varepsilon } 
	\quad {\rm weakly~in~} L^{2}(0, T; H),
	\\
	\label{4.28b}
	\varepsilon \Delta u_{\varepsilon, \lambda_{k}} \to \varepsilon \Delta u_{\varepsilon }
	\quad {\rm weakly~in~} L^{2}(0, T; H)
\end{gather}
as $k\nearrow +\infty $. From \eqref{4.27} and well-known compactness results (see, e.g., \cite{Sim87}), we obtain 
\begin{equation}
	\label{4.33}
	u_{\varepsilon, \lambda_{k}}\to u_{\varepsilon } 
	\quad {\rm strongly~in~} C \bigl( [0, T]; H \bigr)
\end{equation}
as $k\nearrow +\infty $. 
From \eqref{4.33}, we deduce that $u_{\varepsilon}(0)=u_{0\varepsilon }$ a.e.\ in $\Omega $. 
Moreover, from \eqref{4.33} and the {L}ipschitz continuity of $\pi_{\varepsilon }$, we deduce that 
\begin{equation}
	\label{pi}
	\pi_{\varepsilon }(u_{\varepsilon, \lambda_k })
	\to \pi_{\varepsilon }(u_{\varepsilon }) 
	\quad {\rm strongly~in~} C \bigl( [0, T]; H \bigr)
\end{equation}
as $k\nearrow +\infty $. 
Also, applying \eqref{4.30}, \eqref{4.33}, and the monotonicity of $\beta $, 
we obtain $\xi_{\varepsilon}\in \beta (u_{\varepsilon })$ a.e.\ in $Q$. 
Finally, given the level of approximation associated with the weak formulations of \eqref{9}--\eqref{10}, 
taking the limit as $k\nearrow +\infty $, and using \eqref{4.27}--\eqref{pi}, 
we find that $(u_\varepsilon, \mu_\varepsilon, \xi_\varepsilon )$ satisfies \eqref{2.13}--\eqref{2.14}. 
\hfill $\Box$ 

\subsection{Passage to the limit $\varepsilon \searrow 0$}

\paragraph{Proof of Theorem 2.1.}
\indent
From weakly and strongly convergence, \eqref{4.27}--\eqref{4.28b}, the same kind of uniform estimates in Lemmas~3.1 and 3.2 holds for $(u_{\varepsilon }, \mu_{\varepsilon }, \xi_{\varepsilon })$, that is, 
\begin{gather*}
	\int_{0}^{t}
	\bigl|u'_{\varepsilon}(s)
	\bigr|_{V^{*}}^{2}ds
	+\varepsilon \bigl| u_{\varepsilon}(t) \bigr|_{V}^{2}
	+ \frac{c_{1}}{4} \bigl| u_{\varepsilon}(t) \bigr |_{H}^{2} \leq M_{1},\\
	\int_{0}^{t}\bigl| \mu_{\varepsilon}(s) \bigr|_{V}^{2}ds\leq M_{1},\\
	\int_{0}^{t} \bigl|\xi_{\varepsilon}(s) \bigr|_{H}^2ds \leq M_{2}, \\
	\int_{0}^{t}
	\bigl|\varepsilon \Delta u_{\varepsilon}(s)
	\bigr|_{H}^{2}ds \leq M_{3}
\end{gather*}
for all $t \in [0,T]$. 
Hence, there exists a subsequence $\{\varepsilon_{k}\}_{k\in \mathbb{N}}$ with $\varepsilon _{k}\searrow 0$ as $k\nearrow +\infty $ and some limits functions $u\in H^{1}(0, T; V^{*})\cap L^{\infty }(0, T; H)$, $\mu \in L^{2}(0, T; V)$, $\xi \in L^{2}(0, T; H)$ such that 
\begin{gather}
	\label{4.1}
	u_{\varepsilon_{k}}\to u \quad {\rm weakly~star~in~} H^{1}(0, T; V^{*})\cap L^{\infty }(0, T; H), 
	\\ 
	\label{4.2}
	\varepsilon_{k}u_{\varepsilon_{k}}\to 0 \quad {\rm strongly~in~} L^{\infty }(0, T; V), 
	\\
	\label{4.3}
	\mu_{\varepsilon_{k}}\to \mu \quad {\rm weakly~in~} L^{2}(0, T; V), 
	\\
	\label{4.4}
	\xi_{\varepsilon_{k}}\to \xi \quad {\rm weakly~in~} L^{2}(0, T; H)
\end{gather}
as $k\nearrow + \infty $. 
From \eqref{4.1} and the well-known {A}scoli--{A}rzel\`a theorem (see, e.g., \cite{Sim87}), we obtain 
\begin{equation}
	\label{4.5}
	u_{\varepsilon_{k}} \to u \quad {\rm strongly~in~} C \bigl( [0, T]; V^{*} \bigr). 
\end{equation}
Moreover, by virtue of \eqref{4.2},
\begin{equation}
	\varepsilon_{k}\Delta u_{\varepsilon_{k}}
	\to 0 \quad {\rm weakly~in~} L^{2}(0, T; H)
	\label{lap}
\end{equation}
as $k\nearrow +\infty $. 
Aided by assumption {\rm (A3)}, we see that 
\begin{equation*}
	\bigl|\pi_{\varepsilon_{k}}(u_{\varepsilon_{k}})\bigr|
	\leq |\pi_{\varepsilon_{k} }'|_{L^{\infty }(\mathbb{R})}
	|u_{\varepsilon_{k}}|
	+ \bigl| u_{\varepsilon _{k}}(0) \bigr|
	\leq c_{3}\sigma (\varepsilon _{k})\bigl( 1+|u_{\varepsilon_{k}}|\bigr)
	 \quad {\rm a.e.\ in~}Q. 
\end{equation*}
Therefore, we obtain 
\begin{equation}
	\label{4.6}
	\pi_{\varepsilon_{k}}(u_{\varepsilon_{k}})
	\to 0 \quad {\rm strongly~in~} L^{\infty }(0, T; H)
\end{equation}
as $k\nearrow +\infty $. 
Now, from \eqref{2.13}, we have 
\begin{equation*}
	\bigl\langle u_{\varepsilon }'(t), z\bigr\rangle_{V^{*}, V}
	+\bigl\langle F z,\mu_{\varepsilon }(t) \bigr\rangle_{V^*,V} 
	= 0 \quad {\rm for~all~}z\in V,
\end{equation*}
for a.a.\ $t\in (0, T)$. 
Therefore, using \eqref{4.1}, \eqref{4.3} and \eqref{4.5}, 
we obtain \eqref{2.19} by passing to the limit in the above, that is, \eqref{2.19'}. 
On the other hand, in \eqref{2.14}, 
using \eqref{4.3}, \eqref{4.4}, \eqref{lap}, \eqref{4.6} and performing a comparison, 
we obtain $\mu =\xi -f$ in $L^2(0,T;H)$. 
Moreover, by comparison, this gives us the additional regularity 
$\xi =\mu +f \in L^2(0,T;V)$ because $\mu, f \in L^2(0,T;V)$. 
We now have $u\in C([0, T]; V^{*})\cap L^{\infty }(0, T; H)$; 
the function $u$ is thus weakly continuous from $[0, T]$ to $H$, that is, \eqref{ic} holds. 
Finally, it remains to show that $\xi \in \beta (u)$ a.e.\ in $Q$. 
For this purpose, we show that 
\begin{equation}
	\label{4.7}
	\limsup_{k \to +\infty }\int_{0}^{T}
	\bigl(\xi_{\varepsilon_{k}}(t), u_{\varepsilon_{k}}(t)\bigr)_{H}dt
	\leq \int_{0}^{T} \bigl(\xi (t), u(t) \bigr)_{H}dt. 
\end{equation}
To this aim, testing \eqref{2.14} for $u_{\varepsilon_{k}}(t)$ and integrating over $[0, T]$, we deduce from the boundary condition that 
\begin{align*}
	\int_{0}^{T} \bigl(\xi_{\varepsilon_{k}}(t), 
	u_{\varepsilon_{k}}(t)\bigr)_{H}dt
	& = \int_{0}^{T} 
	\bigl(\mu_{\varepsilon_{k}}(t)+f(t), u_{\varepsilon_{k}}(t)
	\bigr)_{H}dt
	+\varepsilon_{k}
	\int_{0}^{T} 
	\bigl(\Delta u_{\varepsilon_{k}}(t), u_{\varepsilon_{k}}(t)\bigr)_{H}dt \\
	& \quad {} - \int_{0}^{T} \bigl( 
	\pi_{\varepsilon_{k}}\bigl( 
	u_{\varepsilon_{k}}(t) 
	\bigr), u_{\varepsilon_{k}} ( t ) 
	\bigr)_{H}dt \\
	& = \int_{0}^{T} 
	\bigl\langle u_{\varepsilon_{k}}(t), \mu_{\varepsilon_{k}}(t)+f(t) \bigr \rangle_{V^*,V} dt 
	-\varepsilon_{k}
	\int_{0}^{T}\bigl|\nabla u_{\varepsilon_{k}}(t)\bigr|_{H^{d}}^{2}dt \\
	& \quad {} 
	- 
	\varepsilon_{k} \kappa  
	\int_{0}^{T} \bigl| u_{\varepsilon_{k}}(t) \bigr|_{H_\Gamma }^2 dt 
	-\int_{0}^{T} \bigl(\pi_{\varepsilon_{k}} \bigl( 
	u_{\varepsilon_{k}}(t)
	\bigr), u_{\varepsilon_{k}}(t)\bigr)_{H}dt \\
	& \leq 
	\int_{0}^{T}\bigl\langle u_{\varepsilon_{k}}(t), 
	\mu_{\varepsilon_{k}}(t)+f(t)
	\bigr\rangle_{V^*,V} dt 
	-\int_{0}^{T}\bigl(\pi_{\varepsilon_{k}}
	\bigl( u_{\varepsilon_{k}}(t) \bigr), u_{\varepsilon_{k}}(t)\bigr)_{H}dt. 
\end{align*}
From \eqref{4.3} and \eqref{4.5}, we see that 
\begin{align*}
	\lim_{k\to \infty }\int_{0}^{T}
	\bigl\langle u_{\varepsilon_{k}}(t), \mu_{\varepsilon_{k}}(t) + f(t)\bigr\rangle_{V^{*}, V}dt
	&= \int_{0}^{T} \bigl\langle u(t), \mu (t)+f(t)\bigr\rangle_{V^{*}, V}dt \\
	&= \int_{0}^{T} \bigl(u(t), \mu (t)+f(t)\bigr)_{H}dt \\
	&= \int_{0}^{T} \bigl(\xi (t), u(t)\bigr)_{H}dt. 
\end{align*}
Moreover, from \eqref{4.1} and \eqref{4.6}, we have 
\begin{equation*}
	\lim_{k\to \infty }\int_{0}^{T}
	\bigl(\pi_{\varepsilon_{k}}
	\bigl( 
	u_{\varepsilon_{k}}(t) 
	\bigr), u_{\varepsilon_{k}}(t)\bigr)_{H}dt
	=0. 
\end{equation*}
Therefore, \eqref{4.7} holds. 
Applying the closedness theorem with respect to weakly-weakly convergence (see, e.g.,\ \cite[Lemma~2.3]{Bar10}), we obtain the fact that $\xi \in \beta (u)$ a.e.\ in $Q$. 
Thus, the proof of existence stated in Theorem~2.1 is complete. 

Next, we show that the component $u$ of the weak solution of (P) is unique. 
Now, for $i=1, 2$, let $(u^{(i)}, \xi ^{(i)})$ be weak solutions of (P), respectively. 
Then, from \eqref{2.19}, we have 
\begin{align*}
	& \bigl
	\langle u_{1}'(t)-u_{2}'(t), 
	z\bigr\rangle_{V^{*}, V}
	+\int_{\Omega }
	\nabla \bigl(\mu_{1}(t)-\mu_{2}(t)\bigr)\cdot \nabla zdx \\
	& \quad {} + \kappa \int_{\Gamma }
	\bigl( \mu_{1}(t)-\mu_{2}(t) \bigr) z d\Gamma =0 \quad {\rm for~all~}z\in V. 
\end{align*}
Here, putting $z:=F^{-1}(u_{1}(t)-u_{2}(t))$ in $V$ at time $t\in (0, T)$, using the monotonicity of $\beta $, and $\mu_{1}-\mu_{2}=\xi_{1}-\xi_{2}$ a.e.\ in $Q$, we infer that 
\begin{align*}
	& \int_{\Omega }
	\nabla \bigl(\mu_{1}(t)-\mu_{2}(t)\bigr)\cdot \nabla F^{-1}\bigl( u_{1}(t)-u_{2}(t) \bigr)dx
	+ \kappa \int_{\Gamma }
	\bigl( \mu_{1}(t)-\mu_{2}(t) \bigr) F^{-1}\bigl( u_{1}(t)-u_{2}(t) \bigr) d\Gamma \\
	& \quad = \bigl\langle FF^{-1} \bigl( u_{1}(t)-u_{2}(t)\bigr), \xi_{1}(t)-\xi_{2}(t)\bigr\rangle_{V^{*}, V} \\
	& \quad = \bigl\langle u_{1}(t)-u_{2}(t), \xi_{1}(t)-\xi_{2}(t) \bigr \rangle_{V^*,V} \\
	& \quad = \bigl(\xi_{1}(t)-\xi_{2}(t),u_{1}(t)-u_{2}(t) \bigr)_H \\
	& \quad \geq 0. 
\end{align*}
Therefore, for all $t\in [0, T]$, we deduce that 
\begin{equation*}
	\frac{1}{2}
	\bigl|u_{1}(t)-u_{2}(t)\bigr|_{V^{*}}^{2}
	+\int_{0}^{t}\bigl(\xi_{1}(t)-\xi_{2}(t),u_{1}(t)-u_{2}(t) \bigr)_Hdt \le 0, 
\end{equation*}
which implies that the component $u$ is unique. 
Thus, we obtain our result. 
\hfill $\Box $

\paragraph{Remark 4.1.} {\rm (i)} 
Allowing the maximal monotone graph $\beta $ to be multi-valued, the component $\xi $ is therefore not uniquely determined. 
However, if $\beta $ is single-valued, then $\xi =\beta (u)$ is unique. \\
{\rm (ii)} The argument for the proof of Theorem~2.1 is essentially the same as in \cite{CF16}. 
Also, we can obtain the error estimates by some reinforcement. 
Indeed, if we add the assumption in \eqref{2.8} of {\rm (A3)} by $\sigma (\varepsilon ):=\varepsilon ^{1/2}$. 
Moreover, if we add the following assumption in {\rm (A5)}, then there exists $c_{4}>0$ such that 
\begin{equation}
	|u_{0\varepsilon }-u_{0}|_{V^{*}}
	\leq c_{4}\varepsilon ^{\frac{1}{4}} \quad {\rm for~all~}\varepsilon \in (0, 1]. 
	\label{ee1}
\end{equation}
Then we obtain error estimate
\begin{equation}
	|u_{\varepsilon }
	-u|_{C([0, T]; V^{*})}^{2}
	+\int_{0}^{T}\bigl(\xi_{\varepsilon }(t)-\xi (t), 
	u_{\varepsilon }(t)-u(t)\bigr)_{H}dt
	\leq C^{*}
	\varepsilon ^{\frac{1}{3}} 
	\label{ee2}
\end{equation}
for all $\varepsilon \in (0, \bar{\varepsilon }]$, where $C^*$ is a positive constant depending only on the data (see, \cite[Theorem~5.1]{CF16}). 

\section{Improvement of the results}
\indent
As mentioned in Introduction, for all of the results from the ``{H}ilbert space approach'' to nonlinear diffusion equation, the growth condition {\rm (A2)} for $\beta $ is a very important assumption. Nevertheless, it is too restricted from the perspective of applications. For this reason, improving the growth condition was studied in \cite{DK99} using the ``lower semicontinuous convex extension''. See also a different approach given in \cite{Aka09}. 
In this section, we consider the improved growth condition {\rm (A2)} for $\beta $.
Hereafter, we assume {\rm (A1)}, and {\rm (A3)} with $\sigma (\varepsilon )=\varepsilon ^{1/2}$, {\rm (A5)}, and {\rm (A6)}; that is, assumption {\rm (A2)} is avoided.

\paragraph{Proof of Theorem 2.2.}
Recall \eqref{10} in the following form
\begin{equation*}
	u_{\varepsilon, \lambda }'(s)+F \mu_{\varepsilon, \lambda}(s)=0 
	\quad {\rm in~} V^*,
\end{equation*}
for a.a.\ $s \in (0,T)$. 
Multiplying $u_{\varepsilon, \lambda }(s)$ by the above, we see that 
\begin{equation}
	\label{6.1}
	\bigl\langle u_{\varepsilon, \lambda }'(s), u_{\varepsilon, \lambda }(s)
	\bigr\rangle_{V^{*}, V}
	+\int_{\Omega }\nabla \mu_{\varepsilon, \lambda }(s)\cdot 
	\nabla u_{\varepsilon, \lambda }(s)dx
	+\kappa \int_{\Gamma }\mu_{\varepsilon, \lambda }(s)u_{\varepsilon, \lambda }(s)d\Gamma =0
\end{equation}
for a.a.\ $s\in (0, T)$. 
Moreover, testing \eqref{9} with $-\Delta u_{\varepsilon, \lambda }(s)$ 
and using boundary conditions and assumption {\rm (A6)}, we deduce that 
\begin{align}
	&\int_{\Omega }\nabla \mu_{\varepsilon, \lambda }(s)
	\cdot \nabla u_{\varepsilon, \lambda }(s)ds 
	+\kappa \int_{\Gamma }\mu_{\varepsilon, \lambda }(s)u_{\varepsilon, \lambda }(s)d\Gamma 
	\nonumber \\
	&=\frac{\lambda }{2}\frac{d}{ds}
	\bigl|\nabla u_{\varepsilon, \lambda }(s)\bigr|_{H^{d}}^{2}
	+\frac{\kappa \lambda }{2}\frac{d}{ds}\bigl|u_{\varepsilon, \lambda }(s)\bigr|_{H_{\Gamma }}^{2}
	+\varepsilon \bigl|\Delta u_{\varepsilon, \lambda }(s)\bigr|_{H}^2
	+\int_{\Omega }\beta '_{\lambda }
	\bigl( u_{\varepsilon, \lambda }(s) \bigr)
	\bigl|
	\nabla u_{\varepsilon, \lambda }(s)
	\bigl|^{2}dx \nonumber \\
	&\quad 
	{}
	+ \kappa 
	\int_{\Gamma }
	u_{\varepsilon, \lambda }(s)
	\beta_{\lambda } \bigl( 
	u_{\varepsilon, \lambda }(s) \bigr) 
	d\Gamma 
	- \bigl(\pi_{\varepsilon } \bigl( 
	u_{\varepsilon, \lambda }(s) \bigr), 
	\Delta u_{\varepsilon, \lambda }(s)\bigr)_{H}
	+\bigl(\Delta f(s), u_{\varepsilon, \lambda }(s)\bigr)_{H} \label{6.2}
\end{align}
for a.a.\ $s\in (0, T)$. Here, we used
\begin{align*}
	\bigl( 
	f(s), \Delta u_{\varepsilon, \lambda }(s)
	\bigr)_{H}
	&=-\int_{\Omega }\nabla f(s)\cdot \nabla u_{\varepsilon, \lambda }(s)dx
	+\int_{\Gamma }f(s) \partial_{\boldsymbol{\nu }} u_{\varepsilon, \lambda }(s)d\Gamma \\
	&=-\int_{\Omega }\nabla f(s)\cdot \nabla u_{\varepsilon, \lambda }(s)dx
	-\kappa \int_{\Gamma } f(s)u_{\varepsilon, \lambda }(s)d\Gamma \\
	&=-\int_{\Omega }\nabla f(s)\cdot \nabla u_{\varepsilon, \lambda }(s)dx
	+\int_{\Gamma }\partial_{\boldsymbol{\nu }}f(s)u_{\varepsilon, \lambda }(s)d\Gamma \\
	&=\bigl(\Delta f(s), u_{\varepsilon, \lambda }(s)\bigr)_{H}
\end{align*}
for a.a.\ $s\in (0, T)$. 
Here, using {\rm (A3)} and the {Y}oung inequality, we infer that 
\begin{align}
	\bigl( \pi_{\varepsilon } \bigl( u_{\varepsilon, \lambda }(s) \bigr), 
	\Delta u_{\varepsilon, \lambda }(s)\bigr)_{H}
	&\leq \int_{\Omega }
	\left(
	|\pi_{\varepsilon }'|_{L^{\infty }(\mathbb{R})}
	\bigl| u_{\varepsilon, \lambda }(s) \bigr|
	+\bigl|\pi_{\varepsilon }(0) \bigr| 
	\right)
	\bigl|
	\Delta u_{\varepsilon, \lambda }(s)
	\bigr|dx \nonumber \\
	&\leq \int_{\Omega }c_{3}
	\varepsilon ^{\frac{1}{2}}
	\left(1+\bigl|u_{\varepsilon, \lambda }(s)\bigr| \right)
	\bigl|\Delta u_{\varepsilon, \lambda }(s)\bigr|dx \nonumber \\
	&\leq \frac{\varepsilon }{2}
	\bigl| \Delta u_{\varepsilon, \lambda }(s) \bigr|_{H}^{2}
	+c_{3}^{2}
	\left( |\Omega |+ \bigl|u_{\varepsilon, \lambda }(s)
	\bigr|_{H}^{2} \right) \label{6.3}
\end{align}
for a.a.\ $s\in (0, T)$. 
Then, combining (\ref{6.1})--(\ref{6.3}) 
and integrating the resultant with respect to $s$ over interval $[0, t]$, we obtain 
\begin{align}
	& \frac{1}{2}\bigl|u_{\varepsilon, \lambda }(t)
	\bigr|_{H}^{2}
	+\frac{\lambda }{2}\bigl|\nabla u_{\varepsilon, \lambda }(t)\bigr|_{H^{d}}^{2}
	+\frac{\kappa \lambda }{2}\bigl|u_{\varepsilon, \lambda }(t)\bigr|_{H_{\Gamma }}^{2}
	+\frac{\varepsilon }{2}
	\int_{0}^{t}\bigl|\Delta u_{\varepsilon, \lambda }(s)\bigr|_{H}^2ds \nonumber \\
	& \quad \leq \frac{1}{2}|u_{0\varepsilon }|_{H}^{2}
	+\frac{\lambda }{2}|\nabla u_{0\varepsilon }|_{H^{d}}^{2}
	+\frac{\kappa \lambda }{2}
	|u_{0\varepsilon }|_{H_{\Gamma }}^{2} \nonumber \\
	& \quad \quad {} +c_{3}^{2}
	\int_{0}^{t}
	\left(
	|\Omega | +
	\bigl|u_{\varepsilon, \lambda }(s)\bigr|_{H}^{2}
	\right)ds
	+\frac{1}{2}\int_{0}^{t}
	\bigl| 
	\Delta f(s) 
	\bigr|_{H}^{2}ds
	+\frac{1}{2}
	\int_{0}^{t}\bigl|u_{\varepsilon, \lambda }(s)\bigr|_{H}^{2}
	ds \label{6.5}
\end{align}
for all $t\in [0, T]$. Now, note that we can take $\lambda \leq \varepsilon $, because $\lambda $ tends to $0$ with fixed $\varepsilon $. Also, from \eqref{2.9} of assumption {\rm (A5)}, we have 
\begin{equation}
	\frac{1}{2}|u_{0\varepsilon }|_{H}^{2}
	+\frac{\lambda }{2}|\nabla u_{0\varepsilon }|_{H^{d}}^{2}
	+\frac{\kappa \lambda }{2}|u_{0\varepsilon }|_{H_{\Gamma }}^{2}\leq \frac{3}{2}(1+\kappa)c_{4}.
	\label{6.5b}
\end{equation}
Then, using \eqref{6.5}, \eqref{6.5b}, and the {G}ronwall inequality, it follows that 
\begin{align*}
	& 
	\bigl|u_{\varepsilon, \lambda }(t)\bigr|_{H}^{2}
	+\lambda \bigl|\nabla u_{\varepsilon, \lambda }(t)\bigr|_{H^{d}}^{2}
	+\kappa \lambda \bigl|u_{\varepsilon, \lambda }(t) \bigr|_{H_{\Gamma }}^{2} \nonumber \\
	& \quad \leq \left\{ 
	3(1+\kappa )c_{4}+ \int_{0}^{T} \bigl| \Delta f(s) \bigr|_H^{2}ds
	+2c_{3}^{2}|\Omega |T \right\} \exp \bigl\{ (2c_{3}^{2}+1)T \bigr\}=:M_4 
	\label{6.6}
\end{align*}
for all $t\in [0, T]$. Moreover, 
\begin{equation*}
	\varepsilon 
	\int_{0}^{t}\bigl|\Delta u_{\varepsilon, \lambda }(s)\bigr|_{H}^2ds 
	\leq M_4(1+2c_3^2 T+ T)
\end{equation*}
for all $t \in [0,T]$. 
Hence, according to Lemma~3.1, using \eqref{3.6} and \eqref{s6} without \eqref{3.7}, we obtain 
\begin{align*}
	&\frac{1}{2}\int_{0}^{t}
	\bigl|u_{\varepsilon, \lambda }'(s)
	\bigr|_{V^{*}}^{2}ds
	+\lambda \int_{0}^{t}
	\bigl|u_{\varepsilon, \lambda }'(s)\bigr|_{H}^{2}ds
	+\frac{\varepsilon }{2}\bigl|u_{\varepsilon, \lambda }(t)\bigr|_{V}^{2}
	+\bigl|
	\widehat{\beta }_{\lambda }
	\bigl( 
	u_{\varepsilon, \lambda }(t)
	\bigr)
	\bigr|_{L^{1}(\Omega )} \nonumber \\
	& \quad \leq 
	\frac{3}{2}c_{4}+\frac{c_1}{8}(1+c_4)+\frac{c_1}{8}(1+M_4)
	+\frac{1}{2}|f|_{L^{2}(0, T; V)}^{2} 
\end{align*}
for all $t\in [0, T]$. 
Therefore, we obtain the same kind of uniform estimates \eqref{m1} and \eqref{m2} in Lemma~3.1 without the growth condition {\rm (A2)}. 
The rest of the proof of Theorem 2.2 is the same as that of Theorem~2.1.
$\hfill $ $\Box $ 

\paragraph{Remark 5.1.} 
As stated in Remark~4.1, we also can improve the error estimate, by assuming {\rm (A1)}, {\rm (A3)} with $\sigma (\varepsilon ):=\varepsilon ^{1/2}$, {\rm (A5)} and \eqref{ee1}.
Then, \eqref{ee2} is improved by stating it in the form 
\begin{equation*}
	|u_{\varepsilon }
	-u|_{C([0, T]; V^{*})}^{2}
	+\int_{0}^{T}\bigl(\xi_{\varepsilon }(t)-\xi (t), 
	u_{\varepsilon }(t)-u(t)\bigr)_{H}dt
	\leq C^{*}
	\varepsilon ^{\frac{1}{2}} 
\end{equation*}
for all $\varepsilon \in (0, \bar{\varepsilon }]$ (see, \cite[Theorem~6.1]{CF16}). 

\section{Asymptotic limits to solutions of {N}eumann problem}

In this section, we establish the order of convergence between the {R}obin problem and the {N}eumann problem related to {\rm (P)}. 

We rewrite our problem for {\rm (P)}, hereafter denoted {\rm (P)}$_{\rm R}$, as follows:
\begin{gather*}
	{\rm (P)}_{\rm R}
	\quad 
	\begin{cases}
	\displaystyle \frac{\partial u_{\kappa }}{\partial t}
	-\Delta \xi_{\kappa }=g_{\kappa }, \quad \xi_{\kappa }\in \beta (u_{\kappa }) 
	\quad {\rm in~}Q, \\
	\partial_{\boldsymbol{\nu }}\xi_{\kappa }+\kappa \xi_{\kappa }
	=h_{\kappa } \quad {\rm on~}\Sigma, \\
	u_{\kappa }(0)=u_{0\kappa } \quad {\rm in~}\Omega 
	\end{cases}
\end{gather*}
for $\kappa >0$, where $g_{\kappa }: Q\to \mathbb{R}$, $h_{\kappa }: \Sigma \to \mathbb{R}$, and $u_{0\kappa }: \Omega \to \mathbb{R}$ are given data, and $\beta $ is same maximal monotone graph as in {\rm (A1)}. 
Moreover, we recall the previous result \cite{CF16} for problem {\rm (P)} subject to the {N}eumann boundary condition,
\begin{gather*}
	{\rm (P)}_{\rm N} 
	\quad 
	\begin{cases}
	\displaystyle 
	\frac{\partial u}{\partial t}
	-\Delta \xi =g, \quad \xi \in \beta (u) \quad {\rm in~}Q, \\
	\partial_{\boldsymbol{\nu }}\xi 
	=h \quad {\rm on~}\Sigma, \\
	u(0)=u_{0} \quad {\rm in~}\Omega 
	\end{cases}
\end{gather*}
which is denoted {\rm (P)}$_{\rm N}$, where $g: Q\to \mathbb{R}$, $h: \Sigma \to \mathbb{R}$, and $u_{0}: \Omega \to \mathbb{R}$ are given data. 
The well-posedness for {\rm (P)}$_{\rm N}$ has already discussed in \cite[Theorem~2.3]{CF16}, specifically, there exists a pair $(u, \xi )$ such that $u\in H^{1}(0, T;V^{*})\cap L^{\infty }(0, T;H)$, $\xi \in L^{2}(0, T;V)$ with $\xi \in \beta (u)$ a.e.\ in $Q$ and they satisfy
\begin{gather}
	\bigl\langle u'(t), z\bigr\rangle_{V^{*}, V}
	+\int_{\Omega }\nabla \xi (t)\cdot \nabla zdx
	=
	\int_{\Omega }^{} g(t)z dx +\int_{\Gamma }^{} h(t)z d\Gamma \nonumber \\
	 \quad {\rm for~all~}z\in V, \label{wfn}\quad 
	 {\rm for~a.a.\ } t \in (0,T), \\
	u(0)=u_{0} \quad {\rm a.e.~in}~\Omega. \nonumber 
\end{gather}

For each $\kappa>0$, we assume in addition that, 
\begin{enumerate}
\item[(A7)] $g$, $g_{\kappa }\in L^{2}(0, T; H)$, $h$, $h_{\kappa }\in L^{2}(0, T; H_{\Gamma })$ and $u_{0}$, $u_{0\kappa }\in H$ with $\widehat{\beta}({u_0}), \widehat{\beta }(u_{0\kappa}) \in L^1(\Omega )$. Moreover, there exists a positive constant $c_{5}$ such that 
\begin{gather}
	|g_{\kappa }-g|_{L^2(0,T;H)} \leq \kappa c_{5}, 
	\quad 
	|h_{\kappa }-h|_{L^2(0,T;H_{\Gamma })}\leq \kappa c_{5}, 
	\quad 
	|u_{0\kappa }-u_{0}|_{V^{*}}\leq \kappa c_{5}.
	\label{7-1}
\end{gather}
\end{enumerate}

Then, we obtain the order of convergence between the {R}obin problem {\rm (P)}$_{\rm R}$
and the {N}eumann problem {\rm (P)}$_{\rm N}$ as $\kappa \searrow 0$:

\paragraph{Theorem 6.1.}
{\it Under assumptions {\rm (A1)}, {\rm (A2)}, {\rm (A7)} or {\rm (A1)}, {\rm (A7)} with $h=h_\kappa \equiv 0$. 
Let $(u, \xi )$ be the weak solution {\rm (P)}$_{\rm N}$ and $(u_{\kappa }, \xi_{\kappa })$ be the weak solution of {\rm (P)}$_{\rm R}$. 
Then, there exists a positive constant $M^{\star }$, depending only the data, such that
\begin{equation}
	|u_{\kappa }-u|_{C([0, T]; V^{*})}^{2}+2\int_{0}^{T}
	\bigl(\xi_{\kappa }(s)-\xi (s), u_{\kappa }(s)-u(s)
	\bigr)_{H}ds\leq M^{\star }\kappa ^{2} \label{pnpr}
\end{equation}
for all $\kappa >0$. Moreover, if $\beta $ is {L}ipschitz continuous, then it follows that
\begin{equation*}
	\int_{0}^{T}
	\bigl| 
	\xi_{\kappa }(s)-\xi (s) \bigr|_{H}^{2}ds\leq C_{\beta }M^{\star }\kappa ^{2}
\end{equation*}
for all $\kappa >0$, where $C_{\beta }$ is the {L}ipschitz constant for $\beta $. }
\paragraph{Proof}
We take the difference between \eqref{2.19'} for $u_\kappa $
and \eqref{wfn} for $u$ at $t=s$ 
and test it setting $z:=F^{-1}(u_{\kappa }(s)-u(s))$ at time $s\in (0, T)$. Then, we have 
\begin{align*}
	&
	\frac{1}{2}\frac{d}{ds}
	\bigl|u_{\kappa }(s)-u(s)\bigr|_{V^{*}}^{2}
	+\int_{\Omega }^{} 
	\nabla \bigl( 
	\xi_{\kappa }(s)-\xi (s) \bigr) \cdot \nabla F^{-1}\bigl( u_{\kappa }(s)-u(s) \bigr) dx \\
	&\quad {}+\kappa \int_{\Gamma }^{} \xi_{\kappa }(s)F^{-1}\bigl( u_{\kappa }(s)-u(s) \bigr) d\Gamma \\
	& =\int_{\Omega }^{} \bigl( g_{\kappa }(s)-g(s) \bigr) 
	F^{-1} \bigl( u_{\kappa }(s)-u(s) \bigr) dx 
	+ \int_{\Gamma }^{} \bigl(h_{\kappa }(s)-h(s) \bigr) 
	F^{-1}(u_{\kappa }(s)-u(s) \bigr) d\Gamma. 
\end{align*}
Now, note that 
\begin{align*}
	&\int_{\Omega }^{} \nabla \bigl( 
	\xi_{\kappa }(s)-\xi (s) \bigr)
	\cdot \nabla F^{-1}\bigl( u_{\kappa }(s)-u(s) \bigr) dx 
	+\kappa \int_{\Gamma }^{} \xi_{\kappa }(s)F^{-1} \bigl( u_{\kappa }(s)-u(s) \bigr) d\Gamma\\
	& =\bigl\langle F F^{-1}\bigl( u_{\kappa }(s)-u(s) \bigr), \xi_{\kappa }(s)-\xi (s)\bigr\rangle_{V^{*}, V}
	+\kappa \bigl(\xi (s), F^{-1}\bigl(u_{\kappa }(s)-u(s)\bigr)\bigr)_{H_{\Gamma }} \\
	& =\bigl(u_{\kappa }(s)-u(s), \xi_{\kappa }(s)-\xi (s) \bigr)_{H}
	+\kappa \bigl(\xi (s), F^{-1}\bigl(u_{\kappa }(s)-u(s)\bigr) \bigr)_{H_{\Gamma }}. 
\end{align*}
Therefore, using the {Y}oung inequality and the trace theorem 
$|\cdot |_{H_\Gamma }^2\le c_{\rm P}'|\cdot |_V^2$, we deduce that
\begin{align}
	&\frac{1}{2}\frac{d}{ds}\bigl|u_{\kappa }(s)-u(s)\bigr|_{V^{*}}^{2}
	+\bigl(\xi_{\kappa }(s)-\xi (s), u_{\kappa }(s)-u(s)\bigr)_{H} \nonumber \\
	&=-\kappa \bigl(\xi (s), F^{-1} \bigl( u_{\kappa }(s)-u(s) \bigr)\bigr)_{H_{\Gamma }}
	+\bigl(g_{\kappa }(s)-g(s), F^{-1} \bigl(u_{\kappa }(s)-u(s) \bigr)\bigr)_{H}
	 \nonumber \\
	&\quad \quad 
	+\bigl(h_{\kappa }(s)-h(s), F^{-1} \bigl(u_{\kappa }(s)-u(s) \bigr)\bigr)_{H_{\Gamma }}
	 \nonumber \\
	&\leq \frac{3}{2}\kappa ^{2}c_{\rm P}'\bigl|\xi (s)\bigr|_{H_{\Gamma }}^{2}
	+\frac{1}{6c_{\rm P}'}\bigl| F^{-1}\bigl( u_{\kappa }(s)-u(s) \bigr)\bigr|_{H_{\Gamma }}^{2}
	+\frac{3}{2} \bigl| g_{\kappa }(s)-g(s)\bigr|_{H}^{2}\nonumber \\
	& \quad {}
	+\frac{1}{6} \bigl| F^{-1}\bigl( u_{\kappa }(s)-u(s) \bigr)\bigr|_{H}^{2} 
	+\frac{3}{2}c_{\rm P}' \bigl|h_{\kappa }(s)-h(s)\bigr|_{H_{\Gamma }}^{2} 
	+\frac{1}{6c_{\rm P}'} \bigl| F^{-1}\bigl( u_{\kappa }(s)-u(s) \bigr)\bigr|_{H_{\Gamma }}^{2} \nonumber \\
	&\leq \frac{3}{2}\kappa ^{2}c_{\rm P}'^2 
	\bigl| \xi (s) \bigr|_{V}^{2}
	+\frac{3}{2} \bigl| g_{\kappa }(s)-g(s)\bigr|_{H}^{2}
	+\frac{3}{2}c_{\rm P}'\bigl|h_{\kappa }(s)-h(s)\bigr|_{H_{\Gamma }}^{2} 
	+\frac{1}{2} \bigl|u_{\kappa }(s)-u(s)\bigr|_{V^{*}}^{2}
	\label{pnpr1}
\end{align}
for a.a.\ $s\in (0, T)$. 
Here, using \eqref{7-1} of {\rm (A7)} and the {G}ronwall inequality, we infer that 
\begin{align*}
	& \bigl| u_{\kappa }(t)-u(t) \bigr|_{V^{*}}^{2} \nonumber \\
	&\leq \left\{ 
	|u_{0\kappa }-u_{0}|_{V^{*}}^{2}
	+3c_{\rm P}'^2 \kappa ^{2}
	| \xi |_{L^2(0,T;V)}^{2} 
	+ 3 |g_{\kappa }-g|_{L^2(0,T;H)}^{2} 
	+ 3c_{\rm P}' |h_{\kappa }-h|_{L^2(0,T;H_\Gamma)}^{2} 
	\right\}\exp (T) \nonumber \\
	&\leq \kappa ^2 \bigl\{c_{5}^{2}
	+3c_{\rm P}'^2| \xi |_{L^2(0,T;V)}^{2} +3c_5^2 + 3 c_{\rm P}'c_5^2\bigr\}\exp (T)
\end{align*}
for all $t\in [0, T]$.
Moreover, by integrating \eqref{pnpr1} over $[0, T]$ with respect $s$, there exists a positive constant $M^{\star }$, such that the estimates \eqref{pnpr} holds. 
If $\beta $ is {L}ipschitz continuous, we infer from its monotonicity that
\begin{align*}
	\int_{0}^{T}
	\bigl|\xi_{\kappa }(s)-\xi (s)\bigr|_{H}^{2}ds
	&\leq \int_{0}^{T}\! \!
	\int_{\Omega }C_{\beta }
	\bigl|u_{\kappa }(s)-u(s)\bigr| \bigl|\xi_{\kappa }(s)-\xi (s)\bigr|dxds \\
	&=C_{\beta }\int_{0}^{T}\bigl(\xi_{\kappa }(s)-\xi (s), u_{\kappa }(s)-u(s)\bigr)_{H}ds \\
	&\leq C_{\beta }M^{\star }\kappa ^{2}.
\end{align*}
This completes the proof. $\hfill $ $\Box $

\section*{Acknowledgments}
We thank Richard Haase, Ph.D, from Edanz Group (www.edanzediting.com/ac) for editing a draft of this manuscript.
\section*{Appendix}

We use the same notation as in the previous sections. 
To characterize the operator $A$ as in Section~3, we employ the maximal monotone theory related to the subdifferential.
To do so, we define a proper lower semicontinuous convex functional $\varphi : H\to [0, +\infty ]$,
\begin{equation*}
	\label{varp}
	\varphi (z):=\begin{cases}
	\displaystyle \frac{1}{2}\int_{\Omega }|\nabla z|^{2}dx
	+\frac{\kappa }{2}\int_{\Gamma }|z|^{2}d\Gamma 
	& {\rm if} 
	\quad z\in V, \\
	+\infty & {\rm otherwise}. \\
\end{cases} 
\end{equation*}

We now present the representation of subdifferential operator $\partial \varphi$ (see e.g., \cite[pp.62--65, Proposition~2.9]{Bar10}).

\paragraph{Lemma A.}
{\it The subdifferential $\partial \varphi $ on $H$ is characterized by}
$$\partial \varphi (z)=-\Delta z \quad {\it in~}H \quad {\it for~all~} z\in D(\partial \varphi )=W. $$

\paragraph{Proof} 
Let $z\in D(\varphi )=V$. For all $z^{*}\in \partial \varphi (z)$ in $H$. From the definition of the subdifferential $\partial \varphi $, we have 
\begin{equation*}
	(z^{*}, \tilde{z})_{H}
	=(\nabla z, \nabla \tilde{z})_{H^{d}}
	+\kappa (z, \tilde{z})_{H_{\Gamma }}
	\quad {\rm for~all~} \tilde{z}\in V.
\end{equation*}
Now, taking $\tilde{z}\in \mathcal{D}(\Omega )$ as the test function, we deduce  $z^{*}=-\Delta z$ in $\mathcal{D}'(\Omega )$, that is,
$z^* = -\Delta z$ in $H$ by comparison. 
Hence, for all $\tilde{z}\in V$, we have 
\begin{equation*}
	(z^{*}, \tilde{z})_{H}
	=(-\Delta z, \tilde{z})_{H} 
	+ \langle \partial_{\boldsymbol{\nu} } z, \tilde{z} \rangle 
	+ \kappa\int_{\Gamma }^{} z \tilde{z} d\Gamma,
\end{equation*}
specifically, by comparison, $\partial_{\boldsymbol{\nu }} z = -\kappa z$ in $H^{1/2}(\Gamma)$. 
Therefore, from the theory of elliptic regularity for {N}eumann problem, it follows that $z \in H^2(\Omega )$. 
Thus, we infer that $\partial \varphi (z)$ is singleton and is equal to $-\Delta z$ with $D(\partial \varphi )=\{z\in H^{2}(\Omega ) : \partial_{\boldsymbol{\nu} }z+\kappa z=0 \ {\rm a.e.~on~}\Gamma \}$. 
\hfill $\Box $ \\

Finally, we can define $A$ by $\partial \varphi $ and then find the same abstract structure for the evolution equation \eqref{1}, which has also been treated in previous works \cite{CF15, KNP95, Kub12}.

\end{document}